\documentclass[final,1p,times]{elsarticle}

\usepackage{amssymb}
\usepackage{amsthm}
\usepackage{amscd}
\usepackage{bm}
\usepackage{amsmath}
\usepackage{amsfonts}
\usepackage{graphicx}
\newtheorem{theorem}{Theorem}[section]
\newtheorem{remark}[theorem]{Remark}
\newtheorem{corollary}[theorem]{Corollary}

\newtheorem{proposition}[theorem]{Proposition}
\usepackage{mathrsfs}
\usepackage{titletoc}
\usepackage{color}

\newcommand{\f}{\frac}

\newcommand{\be}{\begin{equation}}
\newcommand{\ee}{\end{equation}}
\newcommand{\bea}{\begin{eqnarray}}
\newcommand{\eea}{\end{eqnarray}}
\newcommand{\bna}{\begin{eqnarray*}}
\newcommand{\ena}{\end{eqnarray*}}

\renewcommand{\le}{\left}
\newcommand{\ri}{\right}

\usepackage{geometry}
\journal{***}

\begin{document}

\begin{frontmatter}

\title{Sobolev spaces on locally finite graphs}

\author[qnu]{Mengqiu Shao}
\ead{mqshaomath@qfnu.edu.cn}
\address[qnu]{School of Mathematical Sciences, Qufu Normal University, Shandong, 273165, China}

\author[ruc]{Yunyan Yang}
\ead{yunyanyang@ruc.edu.cn}
\address[ruc]{School of Mathematics, Renmin University of China, Beijing, 100872, China}
\author[bnu]{Liang Zhao\corref{Zhao}}
\ead{liangzhao@bnu.edu.cn}
\address[bnu]{School of Mathematical Sciences, Key Laboratory of Mathematics and Complex Systems of MOE,\\
Beijing Normal University, Beijing, 100875, China}

\cortext[Zhao]{Corresponding author.}

\begin{abstract}
In this paper, we develop the theory of Sobolev spaces on locally finite graphs, including completeness, reflexivity, separability, and Sobolev inequalities. Since there is no exact concept of dimension on graphs, classical methods that work on Euclidean spaces or Riemannian manifolds can not be directly applied to graphs. To overcome this obstacle, we introduce a new linear space composed of vector-valued functions with variable dimensions, which is highly applicable for this issue on graphs and is uncommon when we consider to apply the standard proofs on Euclidean spaces to Sobolev spaces on graphs. The gradients of functions on graphs happen to fit into such a space and we can get the desired properties of various Sobolev spaces along this line. Moreover, we also derive several Sobolev inequalities under certain assumptions on measures or weights of graphs. As fundamental analytical tools, all these results would be extremely useful for partial differential equations on locally finite graphs.
\end{abstract}


\end{frontmatter}

\section{Introduction}
The Sobolev spaces, which concern spaces of functions having specific differentiability properties \cite{Sobolev1, Sobolev2}, turn out to be the proper setting in which to apply the methods of functional analysis and have become the fundamental tools in the theory of partial differential equations. As an intersection of PDEs and graph theory, partial differential equations on graphs appear in a wide range of questions in both pure and applied mathematics, such as differential geometry \cite{Horn,Lai}, optimal transport \cite{Chow}, image processing \cite{Ta, Elmoataz}, and have become a subject of great interest. One important direction to study such discrete differential equations is to consider their variational structures and it was first implemented by Grigor'yan, Lin and Yang in their series of papers \cite{Gri1,Gri2,Gri3}. In particular, they have pointed out that the required Sobolev spaces for kinds of nonlinear differential equations on graphs are pre-compact, which makes it possible to apply the variational methods.

To be specific, let $G=(V,E)$ be a connected graph, where $V$ denotes the vertex set and $E$ denotes the edge set. Let $\mu:V\rightarrow (0,+\infty)$ be a measure. For any function $u:V\rightarrow \mathbb{R}$, the Laplacian
of $u$ at $x\in V$ is defined as
\begin{equation}\label{lap}
	\Delta u(x)=\f{1}{\mu(x)}\sum_{y\sim x}w_{xy}(u(y)-u(x)),
\end{equation}
where $w_{xy}$ is the weight for any edge $xy\in E$ and we always assume that $w_{xy}$is strictly positive, $y\sim x$ means $xy\in E$ or $y$ is adjacent to $x$.
Corresponding to the case of Euclidean space, the norm of the gradient of $u$ at $x\in V$ is defined as
\begin{equation}\label{grd}
	|\nabla u|(x)=\le(\f{1}{2\mu(x)}\sum_{y\sim x}w_{xy}(u(y)-u(x))^2\ri)^{1/2}.
\end{equation}
For any $p>0$, $L^p(V)$ denotes a linear space of all functions $u:V\rightarrow\mathbb{R}$ with finite norms
$$
\|u\|_p=\left(\int_V|u|^pd\mu\right)^{1/p},
$$
where the integral of $|u|^p$ over $V$ is
$$
\int_{V}|u|^p d\mu=\sum_{x\in V}\mu(x)|u(x)|^p.
$$
Similarly, $L^\infty(V)$ is a linear space of all bounded functions with the norm
$$
\|u\|_\infty=\sup_{x\in V}|u(x)|.
$$
To introduce the differentiability properties of functions, we define the Sobolev space $W^{1,p}(V)$ as
$$
W^{1,p}(V)=\le\{u:V\rightarrow\mathbb{R}: \int_V(|\nabla u|^p+|u|^p)d\mu<+\infty\ri\},
$$
which is equipped with the norm
$$\|u\|_{W^{1,p}(V)}=\|u\|_{p}+\|\nabla u\|_{p}$$
or if $1<p<+\infty$,  with the equivalent norm
$(\|u\|_p^p+\|\nabla u\|_p^p)^{1/p}$. In particular, the space $W^{1,2}(V)$ is equipped with the inner product
$$\langle u,v\rangle=\int_V\left(uv+\langle\nabla u,\nabla v\rangle\right) d\mu$$
with the associated norm
$$
\|u\|_{W^{1,2}(V)}=\le(\|u\|_2^2+\|\nabla u\|_2^2\ri)^{1/2}.
$$
Finally, to involve higher order derivatives of functions $u:V\rightarrow\mathbb{R}$, we define for any $m\in\mathbb{N}$,
$$\nabla^mu(x)=
\left\{
\begin{array}{lll}
\Delta^ku(x)&{\rm if}& m=2k,\,k=0,1,2,\cdots\\[1.5ex]
\nabla\Delta^ku(x)&{\rm if}& m=2k+1,\,k=0,1,2,\cdots,
\end{array}\right.$$
where $\Delta$ and $\nabla$ are the Laplacian and gradient operators given in \eqref{lap} and \eqref{grd}, $\Delta^ku(x)=\Delta(\Delta^{k-1}u)(x)$ and $\Delta^0u(x)=u(x)$. The set $W^{m,p}(V)$ composed of functions $u:V\rightarrow\mathbb{R}$ such that $|\nabla^j u|\in L^p(V)$ for all $j=0,1,2,\cdots,m$ is also a linear space with the norm
\begin{equation*}
	\|u\|_{W^{m,p}(V)}=\sum_{j=0}^m\|\nabla^j u\|_p
\end{equation*}
or the equivalent norm $\left(\sum_{j=0}^m\|\nabla^j u\|_p^p\right)^{1/p}$ for $1<p<+\infty$.

When $G=(V,E)$ is a finite graph, there have been some literature about Sobolev spaces $W^{m,p}(V)$ and the spectrum of $p$-Laplacian, for example \cite{Chu, Ost, Chang, HuaWang}, which presented some basic properties of the spaces. In fact, for any $m\in \mathbb{N}$ and $p\in (0,+\infty]$, $W^{m,p}(V)$ on a finite graph are all finite dimensional linear normed spaces and their structures and properties are easy to explore. As a consequence, the linear space $W^{m,p}(V)$ can be identified with $\mathbb{R}^k$, where $k$ equals to the number of vertices in $V$. Since every bounded sequence in $\mathbb{R}^k$ has a convergent subsequence and any two norms of a finite dimensional linear space are equivalent, the authors of \cite{Gri2,Gri3} have pointed out that the Sobolev spaces $W^{m,p}(V)$ on a finite graph $G=(V,E)$ are pre-compact. This important observation naturally leads to variational methods when we consider partial differential equations on finite graphs.

The case that $G=(V,E)$ is infinite is more complicated and subtle. In this paper, we will focus on the case of locally finite graphs. Here a graph $G=(V,E)$ is called locally finite if each vertex has finitely many neighbours. The finite graphs are certainly special locally finite graphs, but for brevity, in the rest of the paper, we always assume that there are infinite vertices in the locally finite graph.

There have been some known facts and applications of Sobolev spaces on a locally finite graph. Most of the existing work on a locally finite graph are related to the special situation $m=1$ and $p=2$, namely the space $W^{1,2}(V)$. On one hand, to perform variational analysis on the nonlinear Laplacian equation, firstly we need to explore this suitable working space. For example, let $W_0^{1,2}(V)$ denote the completion of $C_c(V)$, a set of functions with finite support, under the norm $\|\cdot\|_{W^{1,2}(V)}$. The space $W_0^{1,2}(V)$, which is clearly a Hilbert space and thus a reflexive space, becomes the working space of nonlinear PDEs in \cite{Gri3, Zhang-Zhao, Lin-Yang2}, where the existence of nontrivial or ground state solutions on locally finite graphs were proved. On the other hand, the Dirichlet form on graphs is another important topic directly related to $W^{1,2}(V)$. There are many existed work about the Dirichlet form and its applications in the spectrum of Laplacian, Cheeger constant, Sobolev type inequalities and stochastic completeness on locally finite graphs, where we can find many properties of the space $W^{1,2}(V)$. Readers can refer to \cite{Chu, Badr, HuaLin, HuaLi, Keller} and the references therein.

When we turn to the space $W^{1,p}(V)$ ($p\in (1,+\infty]$) on a locally finite graph, in \cite{Han-Shao}, under the assumptions that the measure $\mu$ on $V$ has a positive lower bound $\mu_0>0$, $w_{xy}=w_{yx}$ for any $xy\in E$ and $D(x)=\sum_{y\sim x}w_{xy}/\mu(x)$ is uniformly bounded in $V$, the authors proved that $W^{1,p}(V)$ is a reflexive Banach space if $2<p<+\infty$. Based on this result, the authors of \cite{Han-Shao} considered a nonlinear equation involving the $p$-Laplace operator $\Delta_p$ instead of the usual Laplace operator $\Delta$, where
$$
\Delta_pu(x)=\f{1}{2\mu(x)}\sum_{y\sim x}w_{xy}\le(|\nabla u|^{p-2}(y)+|\nabla u|^{p-2}(x)\ri)(u(y)-u(x))
$$
comes from the variation of $\|\nabla u\|_p^p$. For details, we refer the readers to (\cite{Han-Shao}, Lemma 2.1). However, for more general measures and weights of locally finite graphs, we know little about the space $W^{1,p}(V)$, let alone the space $W^{m,p}(V)$ for $m\geq 2$.

According to the above, our aim in this paper is to investigate properties of $W^{m,p}(V)$, including completeness, reflexivity, and separability, on an arbitrary connected locally finite graph without any other restrictions. Precisely, we first prove the following theorem.

\begin{theorem}\label{wmp}
  Let $G=(V,E)$ be a connected locally finite graph. For any nonnegative integer $m$, the Sobolev space $W^{m,p}(V)$ holds the following properties:

  $(i)$ it is a Banach space for $1\leq p\leq+\infty$;

  $(ii)$ it is reflexive for $1<p<+\infty$;

  $(iii)$ it is separable for $1\leq p<+\infty$.
\end{theorem}

Although the proof of Theorem \ref{wmp} is essentially based on the theory of functional analysis, the situations on discrete infinite graphs are quite different from those on Euclidean spaces or Riemannian manifolds. Here we briefly illustrate the difference. Suppose that $G=(V,E)$ is locally finite and $u\in W^{1,p}(V)$. Then we know that $|u|^p$ and $|\nabla u|^p$ are integrable on $V$. In view of \eqref{grd}, we formally write
$$
\nabla u(x)=\le(\sqrt{\f{w_{x,y_1}}{2\mu(x)}}(u(y_1)-u(x)),\cdots,\sqrt{\f{w_{x,y_{\ell_x}}}{2\mu(x)}}(u(y_1)-u(x))\ri),
$$
where $\ell_x$ denotes the number of all neighbours of $x$. Thus every $\nabla u(x)$ can be viewed as an $\ell_x$ dimensional vector in $\mathbb{R}^{\ell_x}$. Let $\rho(x)=dist(x,O)$ denotes the distance from $x\in V$ to a fixed vertex $O\in V$. In general, $\ell_x$ may tend to infinity as $\rho(x)$ tends to infinity. Hence $\nabla u$ can not be embedded in any finite dimensional vector space. This makes the classical methods unable to be directly applied. To overcome this obstacle, we first introduce a linear space $\mathscr{F}$ composed of vector-valued functions with variable dimensions and a subspace $\mathscr{V}^p$ of $\mathscr{F}$. With the help of the space $\mathscr{V}^p$, we will finally obtain the desired properties of $W^{1,p}(V)$ and $W^{m,p}(V)$ ($m\geq 2$).

We remark that in \cite{Mug}, Mugnolo presented a similar result (Lemma 3.6) with a sketchy and not very readable proof about the space $W^{1,p}_{\rho, \nu}(V)$ on a locally finite graph. Using our notations, the norm of $u\in W^{1,p}_{\rho,\nu}(V)$ in \cite{Mug} is written as follows
\begin{equation*}
	\|u\|_{W^{1,p}_{\rho,\nu}}=\left(\sum_{x\in V}\mu(x)|u(x)|^p \right)^{1/p}+\left(\frac{1}{2}\sum_{x\in V}\sum_{y\sim x}w_{xy}|u(y)-u(x)|^p\right)^{1/p}.
\end{equation*}
Our norm of $u\in W^{1,p}(V)$ is
\begin{equation*}
	\|u\|_{W^{1,p}}=\left(\sum_{x\in V}\mu(x)|u(x)|^p \right)^{1/p}+\left(\sum_{x\in V}\mu(x)\left(\frac{1}{2\mu(x)}\sum_{y\sim x}w_{xy}|u(y)-u(x)|^2\right)^{p/2}\right)^{1/p}.
\end{equation*}
These two norms are obviously not equivalent for $p\neq 2$. As we have mentioned before, since our norm $\|\cdot\|_{W^{1,p}}$ comes from the variation of $\|\nabla \cdot\|_p^p$, it is more natural and useful when we tend to apply variational methods to PDEs on graphs. At the same time, a measure $\mu$ on vertices without a uniformly strict positive lower bound shall cause some essential difficulty when using the norm $\|\cdot\|_{W^{1,p}}$ \cite{Han-Shao}, which will not occur for $\|\cdot\|_{W^{1,p}_{\rho,\nu}}$ in \cite{Mug}.

Another important Sobolev space $W_0^{m,p}(V)$ is defined as a completion of $C_c(V)$ under the norm $\|\cdot\|_{W^{m,p}(V)}$. Since $W_0^{m,p}(V)$ is a closed subspace of $W^{m,p}(V)$, they share some common properties. In particular, we have

\begin{corollary}\label{w0mp}
  Let $G=(V,E)$ be a connected locally finite graph. For any nonnegative integer $m$, the Sobolev space $W_0^{m,p}(V)$
  holds the following properties:

  $(i)$ it is a Banach space for $1\leq p\leq+\infty$;

  $(ii)$ it is reflexive for $1<p<+\infty$;

  $(iii)$ it is separable for $1\leq p<+\infty$.
\end{corollary}

We can even consider some more general Sobolev spaces on locally finite graphs based on properties of the space $\mathscr{V}^p$. Let $m\geq 1$ be an integer and  $\mathcal{P}=(p_0,p_1,\cdots,p_m)$ is a vector with $p_j\geq 1$, $j=0,\cdots,m$. Define a function space
$$W^{m,\mathcal{P}}(V)=\le\{u:V\rightarrow \mathbb{R}: \sum_{j=0}^m\|\nabla^ju\|_{p_j}<+\infty\ri\}$$
with a norm
$$
\|u\|_{W^{m,\mathcal{P}}(V)}=\sum_{j=0}^m\|\nabla^ju\|_{p_j}.
$$
As an analogue of Theorem \ref{wmp}, we have
\begin{theorem}\label{wmpp} Let $G=(V,E)$ be a connected locally finite graph, $m\geq 1$ be an integer and $\mathcal{P}=(p_0,p_1,\cdots,p_m)$. The Sobolev space $W^{m,\mathcal{P}}(V)$ holds the following properties:

$(i)$ it is a Banach space for $1\leq p_j\leq +\infty$, $j=0,\cdots,m$;

$(ii)$ it is reflexive for $1< p_j< +\infty$, $j=0,\cdots,m$;

$(iii)$ it is separable for $1\leq p_j<+\infty$, $j=0,\cdots,m$.
\end{theorem}

Compared with $W_0^{m,p}(V)$ and $W^{m,p}(V)$, we can define $W_0^{m,\mathcal{P}}(V)$ be a completion of $C_c(V)$ under the norm $\|\cdot\|_{W^{m,\mathcal{P}}(V)}$.
Similar to Corollary \ref{w0mp}, we have the following:
\begin{corollary}\label{w0mpp}
Let $G=(V,E)$ be a connected locally finite graph, $m\geq 1$ be an integer and $\mathcal{P}=(p_0,p_1,\cdots,p_m)$. The Sobolev space $W_0^{m,\mathcal{P}}(V)$ holds the following properties:
	
$(i)$ it is a Banach space for $1\leq p_j\leq +\infty$, $j=0,\cdots,m$;

$(ii)$ it is reflexive for $1< p_j< +\infty$, $j=0,\cdots,m$;

$(iii)$ it is separable for $1\leq p_j<+\infty$, $j=0,\cdots,m$.
\end{corollary}

Finally, we are also concerned with the Sobolev inequalities on locally finite graphs. These kinds of inequalities are certainly very useful when we study the corresponding nonlinear equations on locally finite graphs. Since
$$
\mu(x)|f(x)|^p\leq \int_V|f|^pd\mu\ \ {\rm or}\ \ \mu(x)|\mathbf{f}(x)|^p\leq \int_V|\mathbf{f}|^pd\mu,$$
for any $x\in V$ and $f\in L^p(V)$ or $\mathbf{f}\in\mathscr{V}^p$, if $\mu(x)\geq \mu_0>0$ for all $x\in V$, there clearly holds a simple fact that $L^p(V)$ and $\mathscr{V}^p$ are continuously embedded into $L^\infty(V)$ and $\mathscr{V}^\infty$ respectively. Based on the definitions of the gradient operator and the Laplace operator, there have been some other embedding theorems for $W^{m,p}(V)$ under various hypothesis. For example, one can refer to \cite{Han1} and \cite{Lin-Yang2}.
Now let us present a more general embedding for $W^{m,\mathcal{P}}(V)$.

\begin{theorem}\label{mu0} Let $G=(V,E)$ be a connected locally finite graph and $m\geq 1$ be an integer.
Suppose $\mu(x)\geq \mu_0>0$ for all $x\in V$ and $\mathcal{P}=(p_0,p_1,\cdots,p_m)$ with $1\leq p_j\leq+\infty$ for $j=0,\cdots,m$. Then there exists a constant $C$ depending only on $\mu_0$ and $\mathcal{P}$, such that for any $u\in W^{m,\mathcal{P}}(V)$,
\begin{equation}\label{infty-mp}
	\|\nabla^j u\|_\infty\leq C\|u\|_{W^{m,\mathcal{P}}(V)},\quad\forall j=0,\cdots,m.
\end{equation}
Moreover, for $j\in \{0,\cdots,m\}$, if the corresponding $p_{j}$ satisfies $1\leq p_j<+\infty$, there holds for any $q\in[p_j,+\infty)$,
\begin{equation}\label{q-mp}
	\|\nabla^ju\|_q\leq C\|u\|_{W^{m,\mathcal{P}}(V)},
\end{equation}
where $C$ is a constant depending only on $\mu_0$, $p_j$ and $q$.
\end{theorem}
Under stronger assumptions that $\mu(x)$ has both positive lower bound and upper bound, and $w_{xy}$ has a positive lower bound, there holds an embedding for $W_0^{m,1}(V)$.
\begin{theorem}\label{w0} Let $G=(V,E)$ be a connected locally finite graph and $m\geq 1$ be an integer.
Suppose there exist positive constants $\mu_0$, $\mu_1$ and $w_0$ such that
$0<\mu_0\leq\mu(x)\leq \mu_1$ and $w_{xy}\geq w_0>0$ for all $x\in V$ and $xy\in E$. Then there exists a constant $C$ depending only on
$\mu_0$, $\mu_1$, $w_0$ and $m$, such that for any $u\in  W_0^{m,1}(V)$ and any integer $j\in [0, m/2]$,
\begin{equation}\label{Sob-2}
	\|\Delta^ju\|_\infty\leq C\int_V|\nabla\Delta^ju|d\mu.
\end{equation}
\end{theorem}

The remaining part of this paper is organized as follows. Firstly, we prove several properties of the space $L^p(V)$ in Section \ref{sec2}.
Section \ref{sec3} is the foremost part, in which we study the space $\mathscr{V}^p$. Based on the results of Section \ref{sec3}, in Section \ref{sec4}, we prove the desired results for both the Sobolev spaces $W^{m,p}(V)$ and $W_0^{m,p}(V)$ and the general Sobolev spaces $W^{m,\mathcal{P}}(V)$ and $W_0^{m,\mathcal{P}}(V)$. Some kinds of embeddings of Sobolev spaces are discussed in Section \ref{sec5}. Finally, in Section \ref{sec6}, we present some questions which deserve further consideration.

\section{$L^p(V)$ space}\label{sec2}

We begin our study of Sobolev spaces with $L^p(V)$ space on a connected locally finite graph $G=(V,E)$.  Although some of its properties are known by experts, for the convenience of readers, we give the details of their proofs.

\begin{proposition}\label{Lp}
 For any $1\leq p\leq +\infty$, $L^p(V)$ is a Banach space.
\end{proposition}
{\it Proof.}
{\it Case $1$}. $p=+\infty$.

Given a Cauchy sequence $\{u_n\}$ in $L^\infty(V)$. For any $\epsilon>0$, there exists some integer $N>0$ such that
\begin{equation*}
	\|u_n-u_m\|_\infty=\sup_{x\in V}|u_n(x)-u_m(x)|<\epsilon,\quad\forall m,n\geq N.
\end{equation*}
Since $\{u_n\}$ is bounded in $L^\infty(V)$, we can take a subsequence $\{u_{n_k}\}$ such that $\{u_{n_k}(x)\}$ converges to some $a_x\in \mathbb{R}$ for any $x\in V$. Define a function $u(x)=a_x$. Clearly we have $u\in L^\infty(V)$ and  for any $R>0$ there holds that
$$
\sup_{x\in B_R}|u_n(x)-u(x)|= \lim_{k\rightarrow\infty}\sup_{x\in B_R}|u_n(x)-u_{n_k}(x)|\leq \epsilon,\quad\forall n\geq N,
$$
where $B_R=\{x\in V:\rho(x)<R\}$, and $\rho(x)$ denotes the distance between $x$ and some fixed point $O\in V$.
Taking $R\rightarrow+\infty$ and $\epsilon\rightarrow 0^+$, we can conclude that $\lim_{n\rightarrow\infty}u_n=u$ in $L^\infty(V)$.

{\it Case 2.} $1\leq p<+\infty$.

Suppose $\{u_n\}$ is a Cauchy sequence in $L^p(V)$. Then for any $\epsilon>0$, there exists a positive integer $N$
such that
\begin{equation}\label{cauchy2}
	\|u_n-u_m\|_p^p=\sum_{x\in V}\mu(x)|u_n(x)-u_m(x)|^p<\epsilon,\quad\forall n,m\geq N.
\end{equation}
In particular, there exists an integer $N_1>0$ such that
$$\|u_n-u_{N_1}\|_p^p=\sum_{x\in V}\mu(x)|u_n(x)-u_{N_1}(x)|^p<1,\quad\forall n\geq N_1.$$
It follows that $\{u_n\}$ is bounded in $L^p(V)$. Since $V$ is locally finite and have at most countable vertices, there exist a function $u:V\rightarrow\mathbb{R}$ and a subsequence $\{u_{n_k}\}$ of $\{u_{n}\}$ such that as $k\rightarrow+\infty$,
\begin{equation}\label{sub2}
	u_{n_k}(x)\rightarrow u(x),\quad \forall x\in V.
\end{equation}
For any $R>1$, there obviously holds
$$
\int_{B_R}|u|^pd\mu=\lim_{k\rightarrow\infty}\int_{B_R}|u_{n_k}|^pd\mu\leq\limsup_{n\rightarrow\infty}\|u_n\|_p^p,
$$
which gives $u\in L^p(V)$. On the other hand,
combining \eqref{cauchy2} and \eqref{sub2}, we have for all $R>1$,
$$\int_{B_R}|u_n(x)-u(x)|^pd\mu=\sum_{x\in B_R}\mu(x)|u_n(x)-u(x)|^p\leq\epsilon,\quad\forall n\geq N.$$
By taking $R\rightarrow+\infty$, we conclude that
$$\|u_n-u\|_p^p=\sum_{x\in V}\mu(x)|u_n(x)-u(x)|^p\leq \epsilon,\quad\forall n\geq N.$$
This implies that $\{u_n\}$ converges to $u$ in $L^p(V)$. $\hfill\Box$

\begin{proposition}\label{Separable}
For any
$1\leq p<+\infty$, $L^p(V)$ is a separable space and $C_c(V)$ is dense in $L^p(V)$.
\end{proposition}
{\it Proof.} We first show that $C_c(V)$ is dense in $L^p(V)$ for $1\leq p<+\infty$. For any $u\in L^p(V)$, we can define a sequence of functions
$$u_k(x)=\le\{\begin{array}{lll}
u(x)&{\rm if}& x\in B_k\\[1.5ex]
0&{\rm if}& x\in V\setminus B_k,
\end{array}\ri.$$
where $k\in\mathbb{N}$ and $B_k$ is the ball centred at $O\in V$ with radius $k$. Obviously $u_k\in C_c(V)$ and one can easily check that $\{u_k\}$ converges to $u$ in $L^p(V)$.

To prove the separability of $L^p(V)$, we define $u_{k,j}:V\rightarrow\mathbb{Q}$ to be
$$u_{k,j}(x)=\le\{\begin{array}{lll}
a_{k,j}(x)&{\rm if}&x\in B_k\\[1.5ex]
0&{\rm if}&x\in V\setminus B_k,
\end{array}\ri.$$
where $\mathbb{Q}$ denotes the set of all rational numbers and $a_{k,j}(x)\rightarrow u_k(x)$ as $j\rightarrow\infty$ for any $x\in B_k$.
Clearly $\{u_{k,j}\}$ is also dense in $L^p(V)$ and has a countable cardinality. Consequently $L^p(V)$ is separable.
$\hfill\Box$

\begin{proposition}\label{Reflexive}
 For any
$1<p<+\infty$, $L^p(V)$ is reflexive.
\end{proposition}
{\it Proof.} Lemma 5.5 in \cite{Han-Shao} tells us that $L^p(V)$ is a uniformly convex space. According to Theorem 3.31 in \cite{Brezis}, any uniformly convex Banach space is reflexive and the desired result is proved.$\hfill\Box$\\

\begin{remark}
	(i) We remark that it is easy to check the assumptions $(G_1):\sup_{x\in V}\frac{\sum_{y}w_{xy}}{\mu(x)}<+\infty$ and $(G_2): \inf_{x\in V}\mu(x)>\mu_0>0$ in \cite{Han-Shao} is not used in the proof of Lemma 5.5. Therefore, we can use this conclusion directly in our situation. (ii) Since $L^p(V)=W^{0,p}(V)$, combining Propositions \ref{Lp}-\ref{Reflexive}, we can conclude Theorem \ref{wmpp} for the simplest case $m=0$.
\end{remark}

At the end of this section, we additionally provide the Riesz representation formulas on locally finite graphs, which may have an independent interest.

\begin{proposition}\label{Riesz-representation}
For any $\phi$ in the dual space $(L^p(V))^\ast$ of $L^p(V)$, there exists a unique function $u\in L^q(V)$ such that
$$
\langle\phi,f\rangle=\int_Vufd\mu,\quad\forall f\in L^p(V),
$$
where $1<p<+\infty$ and $1/p+1/q=1$. Moreover, there holds $\|u\|_q=\|\phi\|_{(L^p(V))^\ast}$.
\end{proposition}
{\it Proof.} The argument is almost the same as that of the Euclidean case and one can refer to Theorem 4.11 in \cite{Brezis}. For the convenience of readers, we also present its details here.

Define an operator $T: L^q(V)\rightarrow (L^p(V))^\ast$ by
$$
\langle Tu,f\rangle=\int_Vufd\mu,\quad\forall u\in L^q(V), \forall f\in L^p(V).
$$
H\"older's inequality gives us
\begin{equation}\label{leq}\|Tu\|_{(L^p(V))^\ast}=\sup_{\|f\|_p= 1}\langle Tu,f\rangle\leq \|u\|_q.\end{equation}
On the other hand, for any fixed $u\in L^q(V)$ with $u\not\equiv 0$, let $f=|u|^{q-2}u/\|u\|_q^{q-1}$. Obviously, we have $\|f\|_p=1$ and thus
\begin{equation}\label{geq}
	\|Tu\|_{(L^p(V))^\ast}\geq \langle Tu,f\rangle=\|u\|_q.
\end{equation}
Combining \eqref{leq} and \eqref{geq}, we obtain $\|Tu\|_{(L^p(V))^\ast}=\|u\|_q$ for all $u\in L^q(V)$. Hence $T$ is injective and
continuous.

We now prove that $T$ is also surjective. Let $E=T(L^q(V))\subset (L^p(V))^\ast$. Since $E$ is closed, it is sufficient to prove that $E$ is dense
in $(L^p(V))^\ast$. Take $h\in (L^p(V))^{\ast\ast}$ satisfying $\langle h,Tu\rangle=0$ for all $u\in L^q(V)$. Since $L^p(V)$
is reflexive, we have $h\in L^p(V)$ and
\begin{equation}\label{uh}
	\int_Vuhd\mu=0,\quad\forall u\in L^q(V).
\end{equation}
If $h\nequiv 0$, by taking $u=|h|^{p-2}h$ in \eqref{uh}, we get a contradiction. Thus $h\equiv 0$. As a consequence, $E$ is dense
in $(L^p(V))^\ast$. Therefore $T: L^q(V)\rightarrow (L^p(V))^\ast$ is an isometry, as what is desired. $\hfill\Box$\\

\begin{proposition}\label{L1-L-infty}
For any $\phi\in(L^1(V))^\ast$, there exists a unique
function $u\in L^\infty(V)$ such that
$$\langle\phi,f\rangle=\int_Vufd\mu,\quad\forall f\in L^1(V),$$
and $\|u\|_\infty=\|\phi\|_{(L^1(V))^\ast}$. Moreover, both $L^1(V)$ and $L^\infty(V)$ are not reflexive.
\end{proposition}
{\it Proof.} The proof is just an adaptation of the Euclidean case (\cite{Brezis}, Pages 99-102). Since
we have done such kind of  adaptation in the proof of Proposition \ref{Riesz-representation}, the details are left to readers. $\hfill\Box$

\section{$\mathscr{V}^p$ space}\label{sec3}
At beginning of this section, let us introduce the concept of the vector-valued function space $\mathscr{V}^p$, which turns out to play the key role in the proof of our theorems. For any $x\in V$, we use $\ell_x$ to denote $\sharp\{y\in V: y\sim x\}$, i.e., the number of all distinct neighbours of $x$. Since $V$ is locally finite and connected, $\ell_x$ is a positive integer for any $x\in V$. Define a vector-valued function $\mathbf{f}$ from $V$ to $\cup_{\ell=1}^\infty\mathbb{R}^\ell$ as follows: for any $x\in V$, there exists a unique vector $\mathbf{v}=(v_1,\cdots,v_{\ell_x})\in\mathbb{R}^{\ell_x}\subset \cup_{\ell=1}^\infty\mathbb{R}^\ell$ such that $\mathbf{f}(x)=\mathbf{v}$. Denote the set of all such vector-valued functions by $\mathscr{F}$. We point out that the dimension of the vector $\mathbf{f}(x)$ may change with the change of $x$.

If $1\leq p<+\infty$, we define a specific subset of $\mathscr{F}$ by
$$\mathscr{V}^p=\left\{\mathbf{f}\in \mathscr{F}: \int_V|\mathbf{f}|^pd\mu<+\infty\right\},$$
where $|\mathbf{f}(x)|=|\mathbf{v}|=\sqrt{v_1^2+\cdots+v_{\ell_x}^2}$ is the usual Euclidean norm. For $\mathbf{f},\mathbf{g}\in\mathscr{V}^p$ and $a,b\in\mathbb{R}$, under the algebraic operation $(a\mathbf{f}+b\mathbf{g})(x)=a\mathbf{f}(x)+b\mathbf{g}(x)$ for $x\in V$, $\mathscr{V}^p$ becomes a linear space. Similarly, if $p=+\infty$, the subset $\mathscr{V}^\infty$ of $\mathscr{F}$ is written as
$$\mathscr{V}^\infty=\left\{\mathbf{f}\in\mathscr{F}: \sup_{x\in V}|\mathbf{f}(x)|<+\infty\right\},$$
which is also a linear space.

Moreover, an easy calculation tells us that both the H\"older inequality and the Minkowski inequality hold on locally finite graphs. In particular, for any $\mathbf{f}\in\mathscr{V}^p$ and $\mathbf{g}\in\mathscr{V}^q$, we have
$$\int_V\langle \mathbf{f}(x),\mathbf{g}(x)\rangle d\mu\leq \left(\int_V|\mathbf{f}(x)|^pd\mu\right)^{1/p}\left(\int_V|\mathbf{g}(x)|^qd\mu\right)^{1/q},$$
where $1/p+1/q=1$ and $\langle \cdot,\cdot\rangle$ denotes the inner product in $\mathbb{R}^{\ell_x}$; for any $\mathbf{f}, \mathbf{g}\in {\mathscr{V}^p}$, we have
$$\left(\int_V|\mathbf{f}+\mathbf{g}|^pd\mu\right)^{1/p}\leq \left(\int_V|\mathbf{f}|^pd\mu\right)^{1/p}+\left(\int_V|\mathbf{g}|^pd\mu\right)^{1/p}.$$
Consequently, with the norms
$$\|\mathbf{f}\|_{\mathscr{V}^p}=\|\mathbf{f}\|_p=\left(\int_V|\mathbf{f}|^pd\mu\right)^{1/p}$$
and
$$\|\mathbf{f}\|_{\mathscr{V}^\infty}=\|\mathbf{f}\|_\infty=\sup_{x\in V}|\mathbf{f}(x)|,$$
both $\mathscr{V}^p$ and $\mathscr{V}^\infty$ become linear normed spaces respectively.

Now we start to investigate properties of the space $\mathscr{V}^p$ for $1\leq p\leq +\infty$. Let us first consider the completeness of $\mathscr{V}^p$, namely
\begin{proposition}\label{Sob-v}
For any $1\leq p\leq +\infty$, $\mathscr{V}^p$ is a Banach space.
\end{proposition}
{\it Proof.}
{\it Case $1$}. $p=+\infty$.

Suppose that $\{\mathbf{f}_n\}$ is a Cauchy sequence in $\mathscr{V}^\infty$. For any $\epsilon>0$, there exists some integer $N>0$ such that
\begin{equation*}
	\|\mathbf{f}_n-\mathbf{f}_m\|_\infty=\sup_{x\in V}|\mathbf{f}_n(x)-\mathbf{f}_m(x)|<\epsilon,\quad\forall n,m\geq N.
\end{equation*}
It is easy to obtain that $\{\mathbf{f}_n\}$ is bounded in $\mathscr{V}^\infty$ and we can choose a subsequence $\{\mathbf{f}_{n_k}\}$ such that $\{\mathbf{f}_{n_k}(x)\}$ converges to some vector $\mathbf{f}_x \in \mathbb{R}^{\ell_x}$ for any $x\in V$. Define a vector-valued function $\mathbf{f}(x)=\mathbf{f}_x$. It is obvious that $\mathbf{f}\in \mathscr{V}^\infty$. For any $R>0$, we have
$$\max_{x\in B_R}|\mathbf{f}_n(x)-\mathbf{f}(x)|= \lim_{k\rightarrow\infty}\max_{x\in B_R}|\mathbf{f}_n(x)-\mathbf{f}_{n_k}(x)|\leq \epsilon,\quad\forall n\geq N.$$
Taking $R\rightarrow+\infty$ and $\epsilon\rightarrow 0^+$, we conclude the convergence of $\{\mathbf{f}_n\}$ to $\mathbf{f}$ in $\mathscr{V}^\infty$
as $n\rightarrow\infty$.

{\it Case 2.} $1\leq p<+\infty$.

Suppose that $\{\mathbf{f}_n\}$ is a Cauchy sequence in $\mathscr{V}^p$. Then for any $\epsilon>0$, there exists a positive integer $N$
such that
\begin{equation}\label{cauchy}
	\|\mathbf{f}_n-\mathbf{f}_m\|_p^p=\sum_{x\in V}\mu(x)|\mathbf{f}_n(x)-\mathbf{f}_m(x)|^p<\epsilon,\quad\forall n,m\geq N.
\end{equation}
In particular, there exists an integer $N_1>0$ such that
$$\|\mathbf{f}_n-\mathbf{f}_{N_1}\|_p^p=\sum_{x\in V}\mu(x)|\mathbf{f}_n(x)-\mathbf{f}_{N_1}(x)|^p<1,\quad\forall n\geq N_1.$$
It follows that $\{\mathbf{f}_n\}$ is bounded in $\mathscr{V}^p$. Since $V$ is locally finite and has at most countable vertices, there exist a vector-valued function
$\mathbf{f}\in\mathscr{F}$ and a subsequence $\{\mathbf{f}_{n_k}\}$ such that as $k\rightarrow+\infty$,
\begin{equation}\label{sub}
	\mathbf{f}_{n_k}(x)\rightarrow \mathbf{f}(x)\quad{\rm in}\quad \mathbb{R}^{\ell_x},\quad \forall x\in V.
\end{equation}
On one hand, for all $R>1$, there holds
$$\int_{B_R}|\mathbf{f}|^pd\mu=\lim_{k\rightarrow\infty}\int_{B_R}|\mathbf{f}_{n_k}|^pd\mu\leq\limsup_{n\rightarrow\infty}\|\mathbf{f}_n\|_p^p,
$$
which gives $\mathbf{f}\in \mathscr{V}^p$. On the other hand,
combining \eqref{cauchy} and \eqref{sub}, we have for all $R>1$,
$$\int_{B_R}|\mathbf{f}_n(x)-\mathbf{f}(x)|^pd\mu=\sum_{x\in B_R}\mu(x)|\mathbf{f}_n(x)-\mathbf{f}(x)|^p\leq\epsilon,\quad\forall n\geq N.$$
Taking $R\rightarrow+\infty$, we conclude that
$$\|\mathbf{f}_n-\mathbf{f}\|_p^p=\sum_{x\in V}\mu(x)|\mathbf{f}_n(x)-\mathbf{f}(x)|^p\leq \epsilon,\quad\forall n\geq N,$$
which gives us the convergence of $\{\mathbf{f}_n\}$ to $\mathbf{f}$ in $\mathscr{V}^p$. $\hfill\Box$

\begin{proposition}\label{reflex-v}
For any $1<p<+\infty$, $\mathscr{V}^p$ is a reflexive space.
\end{proposition}
{\it Proof.} Note that $\mathscr{V}^p$ is already a Banach space by Proposition \ref{Sob-v}.  We divide the proof of this proposition into two steps.

{\it Step 1.} $\mathscr{V}^p$ is reflexive for $2\leq p<+\infty$.

Since the function $(x^2+1)^{p/2}-x^p-1$ is increasing with respect to $x\in[0,+\infty)$, we have
\begin{equation}\label{ab}
	a^p+b^p\leq (a^2+b^2)^{p/2},\quad\forall a, b\geq 0.
\end{equation}
For some positive integer $\ell$, take
$$
a=\left|\frac{\mathbf{v}+\mathbf{w}}{2}\right|,\quad b=\left|\frac{\mathbf{v}-\mathbf{w}}{2}\right|,\quad\forall \mathbf{v},\mathbf{w}\in \mathbb{R}^\ell.$$
Then \eqref{ab} gives us
\begin{eqnarray*}
\left|\frac{\mathbf{v}+\mathbf{w}}{2}\right|^p+\left|\frac{\mathbf{v}-\mathbf{w}}{2}\right|^p&\leq& \left(\left|\frac{\mathbf{v}+\mathbf{w}}{2}\right|^2+\left|\frac{\mathbf{v}-\mathbf{w}}{2}\right|^2\right)^{p/2}\\
&=&\left(\frac{|\mathbf{v}|^2}{2}+\frac{|\mathbf{w}|^2}{2}\right)^{p/2}\\
&\leq & \frac{1}{2}(|\mathbf{v}|^p+|\mathbf{w}|^p),
\end{eqnarray*}
where the last inequality follows from the convexity of the function $\phi(x)=x^{p/2}$ when $x\geq 0$ and $p\geq 2$.
This immediately leads to
\begin{equation}\label{Clarkson}
	\left\|\frac{\mathbf{f}+\mathbf{g}}{2}\right\|_p^p+\left\|\frac{\mathbf{f}-\mathbf{g}}{2}\right\|_p^p\leq \frac{1}{2}\left(\|\mathbf{f}\|_p^p+\|\mathbf{g}\|_p^p\right),\quad\forall \mathbf{f},\mathbf{g}\in \mathscr{V}^p.
\end{equation}

Take $\epsilon>0$ and $\mathbf{f},\mathbf{g}\in\mathscr{V}^p$ satisfying
$\|\mathbf{f}\|_p\leq 1$, $\|\mathbf{g}\|_p\leq 1$ and $\|\mathbf{f}-\mathbf{g}\|_p>\epsilon$. In view of \eqref{Clarkson}, we have
$$
\left\|\frac{\mathbf{f}+\mathbf{g}}{2}\right\|_p^p<1-\left(\frac{\epsilon}{2}\right)^p,
$$
which implies that
$$\left\|\frac{\mathbf{f}+\mathbf{g}}{2}\right\|_p<1-\delta,$$
where
$$\delta=1-\left(1-\left(\frac{\epsilon}{2}\right)^p\right)^{1/p}>0.$$
Therefore, $\mathscr{V}^p$ is uniformly convex and thus reflexive by Theorem 3.31 in \cite{Brezis}.

{\it Step 2.} $\mathscr{V}^p$ is reflexive for $1<p\leq 2$.

Suppose that $1<p<+\infty$ and $1/p+1/q=1$. For any fixed $\mathbf{f}\in \mathscr{V}^p$, the mapping $\mathbf{g}\in\mathscr{V}^q\mapsto \int_V\langle \mathbf{f},\mathbf{g}\rangle d\mu$ is
a continuous linear functional on $\mathscr{V}^q$ and thus it gives an element in $(\mathscr{V}^q)^\ast$. Along this way, we obtain an operator $T: \mathscr{V}^p\rightarrow (\mathscr{V}^q)^\ast$ such that
$$
\langle T\mathbf{f},\mathbf{g}\rangle=\int_V\langle \mathbf{f},\mathbf{g}\rangle d\mu,\quad\forall \mathbf{g}\in\mathscr{V}^q.
$$
We claim that
\begin{equation}\label{norm}
	\|T\mathbf{f}\|_{(\mathscr{V}^q)^\ast}=\|\mathbf{f}\|_p,\quad \forall \mathbf{f}\in\mathscr{V}^p.
\end{equation}
If $\mathbf{f}(x)=\mathbf{0}\in \mathbb{R}^{\ell_x}$ for any $x\in V$, the claim is trivial and in the rest of the proof, we assume that at least $\mathbf{f}(x)\neq\mathbf{0}$ at one vertex $x\in V$. It follows from H\"older's inequality that
$$|\langle T\mathbf{f},\mathbf{g}\rangle|\leq \|\mathbf{f}\|_p\|\mathbf{g}\|_q,\quad\forall \mathbf{g}\in\mathscr{V}^q,$$
and thus
\begin{equation}\label{ll}
	\|T\mathbf{f}\|_{(\mathscr{V}^q)^\ast}\leq \|\mathbf{f}\|_p.
\end{equation}
On the other hand, define a function $\mathbf{g}_0$ as
$$
\mathbf{g}_0(x)=\left\{\begin{array}{lll}
|\mathbf{f}(x)|^{p-2}\mathbf{f}(x)&{\rm if}& \mathbf{f}(x)\not=\mathbf{0}\\[1.5ex]
\mathbf{0}&{\rm if}& \mathbf{f}(x)=\mathbf{0}.
\end{array}\right.
$$
Obviously, we have $\mathbf{g}_0\in \mathscr{V}^q$, $\|\mathbf{g}_0\|_q=\|\mathbf{f}\|_p^{p-1}$ and $\langle T\mathbf{f},\mathbf{g}_0\rangle=\|\mathbf{f}\|_p^p$. Hence
$$
\|T\mathbf{f}\|_{(\mathscr{V}^q)^\ast}\geq \frac{\langle T\mathbf{f},\mathbf{g}_0\rangle}{\|\mathbf{g}_0\|_q}=\|\mathbf{f}\|_p.
$$
This together with \eqref{ll} implies that $T:\mathscr{V}^p\rightarrow (\mathscr{V}^q)^\ast$ is an isometry. In particular, since $\mathscr{V}^p$ is a Banach space, $T(\mathscr{V}^p)$ is a closed subspace of $(\mathscr{V}^q)^\ast$.

We now assume $1<p\leq 2$, $2\leq q<+\infty$ and $1/p+1/q=1$. Firstly, $\mathscr{V}^q$ is reflexive by Step 1. Secondly, since the dual space of a reflexive Banach space is also reflexive (\cite{Brezis}, Corollary 3.21), $(\mathscr{V}^q)^\ast$ is reflexive. Finally, since any closed subspace of a reflexive space is reflexive (\cite{Brezis}, Proposition 3.20), $T(\mathscr{V}^p)$ is reflexive, which implies that $\mathscr{V}^p$ is also reflexive.
$\hfill\Box$

\begin{proposition}\label{separable}
For any $1\leq p<+\infty$, $\mathscr{V}^p$ is a separable space.
\end{proposition}
{\it Proof.} For any $\mathbf{f}\in \mathscr{F}$, we denote its support by ${\rm supp} \mathbf{f}=\{x\in V: \mathbf{f}(x)\not=\mathbf{0}\}$ and
$$
\mathbf{C}_c(V)=\{\mathbf{f}\in\mathscr{F}: {\rm supp} \mathbf{f}\,\, {\rm is\,\,finite}\}.
$$
Clearly $\mathbf{C}_c(V)$ is dense in $\mathscr{V}^p$ for any $1\leq p<+\infty$. Namely, for any $\mathbf{f}\in \mathscr{V}^p$, there exists a vector-valued function sequence $\{\mathbf{f}_{k}\}\subset \mathbf{C}_c(V)$ such that $\|\mathbf{f}_{k}-\mathbf{f}\|_{p}\rightarrow 0$ as $k\rightarrow\infty$. For any $k\in \mathbb{N}$, we can define a sequence of vector-valued functions $\{\mathbf{f}_{k,j}: V\rightarrow\cup_{x\in V}\mathbb{Q}^{\ell_x}\}$ as
$$
\mathbf{f}_{k,j}(x)=\left\{\begin{array}{lll}
\mathbf{b}_{k,j}(x)&{\rm if}&x\in B_k\\[1.5ex]
\mathbf{0}&{\rm if}&x\in V\setminus B_k,
\end{array}\right.$$
where $\mathbb{Q}^{\ell_x}=\{\mathbf{v}=(v_{1},v_{2},\cdots,v_{\ell_x}):v_{i}\in\mathbb{Q}, i=1,\cdots,\ell_x\}\subset \mathbb{R}^{\ell_x}$, and $\mathbf{b}_{k,j}(x)\rightarrow \mathbf{f}_k(x)$ as $j\rightarrow\infty$ for any $x\in B_k$.
Clearly $\{\mathbf{f}_{k,j}\}$ is still dense in $\mathscr{V}^p$ and has a countable cardinality. Consequently $\mathscr{V}^p$ is separable.
$\hfill\Box$\\

Similar to the $L^p(V)$ space, there also holds the Riesz representation formulas for the $\mathscr{V}^p$ space.

\begin{proposition}\label{Riesz}
For any $\mathbf{f}^\ast\in(\mathscr{V}^p)^\ast$, there exists a unique vector-valued function $\mathbf{f}\in\mathscr{V}^q$ such that
$$
\langle\mathbf{f}^\ast,\mathbf{g}\rangle=\int_V\langle \mathbf{f},\mathbf{g}\rangle d\mu,\quad \forall \mathbf{g}\in\mathscr{V}^p,
$$
where $1<p<+\infty$ and $1/p+1/q=1$. Moreover, there holds $\|\mathbf{f}\|_{q}=\|\mathbf{f}^\ast\|_{(\mathscr{V}^p)^\ast}$.
\end{proposition}
{\it Proof.} Let $T: \mathscr{V}^p\rightarrow (\mathscr{V}^q)^\ast$ be defined as in Step 2 of the proof of Proposition \ref{reflex-v}.
Since we have already proved \eqref{norm}, we only need to prove that $T$ is surjective, which is equivalent to show that
$T(\mathscr{V}^p)$ is dense in $(\mathscr{V}^q)^\ast$. In fact, for any $\mathbf{h}^{\ast\ast}\in (\mathscr{V}^q)^{\ast\ast}$ and $\mathbf{h}^{\ast\ast}\nequiv \mathbf{0}$, since $\mathscr{V}^q$ is reflexive, we have $\mathbf{h}^{\ast\ast}\in \mathscr{V}^q$. If $\mathbf{h}^{\ast\ast}$ satisfies that $\langle \mathbf{h}^{\ast\ast},T\mathbf{g}\rangle=0$ for all $\mathbf{g}\in \mathscr{V}^p$, there holds
$$
\int_V\langle \mathbf{h}^{\ast\ast},\mathbf{g}\rangle d\mu=0,\quad\forall \mathbf{g}\in \mathscr{V}^p.
$$
By choosing $\mathbf{g}=|\mathbf{h}^{\ast\ast}|^{p-2}\mathbf{h}^{\ast\ast}\in \mathscr{V}^q$, we obtain $\mathbf{h}^{\ast\ast}=\mathbf{0}$. $\hfill\Box$

\begin{proposition}\label{L1-dual}
For any $\mathbf{f}^\ast\in(\mathscr{V}^1)^\ast$, there exists a unique
$\mathbf{f}\in\mathscr{V}^\infty$ such that 
$$\langle\mathbf{f}^\ast,\mathbf{g}\rangle=\int_V\langle \mathbf{f},\mathbf{g}\rangle d\mu,\quad\forall \mathbf{g}\in \mathscr{V}^1,$$
and $\|\mathbf{f}\|_\infty=\|\mathbf{f}^\ast\|_{(\mathscr{V}^1)^\ast}$. Moreover, both $\mathscr{V}^1$ and $\mathscr{V}^\infty$ are not reflexive.
\end{proposition}
{\it Proof.} Similar to the proof of Proposition \ref{L1-L-infty}, we omit the details and refer readers to (\cite{Brezis}, Pages 99-102). $\hfill\Box$

\section{$W^{m,p}(V)$ and $W^{m,\mathcal{P}}(V)$}\label{sec4}

In this section, to use the powerful tool $\mathscr{V}^p$ in the study of $W^{m,p}(V)$ and $W^{m,\mathcal{P}}(V)$, we shall look at the gradient of function on a locally finite graph from the following point of view. For any fixed point $x\in V$, we assume that the set $\{y_1,\cdots,y_{\ell_x}\}$ contains all its neighbours. In view of (\ref{grd}), we may write the gradient of $u$ at $x$ by
$$\nabla u(x)=\le(\sqrt{\f{w_{xy_1}}{2\mu(x)}}(u(y_1)-u(x)),\cdots,\sqrt{\f{w_{xy_{\ell_x}}}{2\mu(x)}}(u(y_{\ell_x})-u(x))\ri).$$
Clearly $\nabla u(x)\in \mathbb{R}^{\ell_x}$ and satisfy
$$|\nabla u|^p(x)=\le(\f{1}{2\mu(x)}\sum_{y\sim x}w_{xy}(u(y)-u(x))^2\ri)^{p/2},$$
where $|\cdot|$ denotes the usual norm in $\mathbb{R}^{\ell_x}$. For higher order derivatives, if $j\in\mathbb{N}^\ast=\mathbb{N}\setminus\{0\}=\{1,2,\cdots\}$ is even, then $\nabla^ju(x)=\Delta^{j/2}u(x)$ is just a scalar function. If $j\in\mathbb{N}^\ast$ is odd, then $\nabla^ju(x)=\nabla\Delta^{\f{j-1}{2}}u(x)\in \mathbb{R}^{\ell_x}$ is a vector-valued function similar to $\nabla u(x)$. With these facts in our mind, we can begin the proof of Theorem \ref{wmp}.\\

{\it Proof of Theorem \ref{wmp}.} We first deal with the case $m=2k$ for some integer $k\geq 1$.

$(i)$ {\it Completeness for $1\leq p\leq +\infty$.}

Suppose that $\{u_n\}$ is a Cauchy sequence in $W^{2k,p}(V)$. Then if $j$ is even and $0\leq j\leq 2k$, $\{\Delta^{j/2}u_n\}$ is a Cauchy sequence in $L^p(V)$. If $j$ is odd and $1\leq j\leq 2k-1$, $\{\nabla\Delta^{\f{j-1}{2}}u_n\}$ is a Cauchy sequence in $\mathscr{V}^p$. Since both $L^p(V)$ and $\mathscr{V}^p$ are Banach spaces for $1\leq p\leq+\infty$ (Propositions \ref{Lp} and \ref{Sob-v}), we have that $u_n$ converges to some $u$ in $L^p(V)$. Moreover, if $j$ is even and $2\leq j\leq 2k$, $\Delta^{\f{j}{2}}u_n$ converges to some $f_j$ in $L^p(V)$, or if $j$ is odd and  $1\leq j\leq 2k-1$, $\nabla\Delta^{\f{j-1}{2}}u_n$ converges to some $\mathbf{g}_j$ in $\mathscr{V}^p$. It follows from $u_n\rightarrow u$ in $L^p(V)$ that $u_n(x)$ converges to $u(x)$ pointwisely in $x\in V$, which immediately leads to $f_j(x)=\Delta^{\f{j}{2}}u(x)$ and $\mathbf{g}_j(x)=\nabla\Delta^{\f{j-1}{2}}u(x)$ for all $x\in V$. Hence $u\in W^{2k,p}(V)$ and $u_n$ converges to $u$ in $W^{2k,p}(V)$.

$(ii)$ {\it Reflexivity for $1<p<+\infty$.}

Let $E=L^p(V)\times \mathscr{V}^p\times\cdots\times L^p(V)\times \mathscr{V}^p\times L^p(V)$ be a product space of $(k+1)$ $L^p(V)$ spaces and $k$ $\mathscr{V}^p$ spaces. Clearly, $E$ is reflexive. Define an operator $\Gamma:W^{2k,p}(V)\rightarrow E$ as $\Gamma u=(u,\nabla u,\cdots,\Delta^ku)$. This operator satisfies $\|\Gamma u\|_E=\sum_{j=0}^{2k}\|\nabla^ju\|_p=\|u\|_{W^{2k,p}(V)}$ and thus $\Gamma:W^{2k,p}(V)\rightarrow E$ is an isometry. Since $W^{2k,p}(V)$ is a Banach space, $\Gamma(W^{2k,p}(V))$ is a closed subspace of $E$. By Proposition 3.20 in \cite{Brezis}, we conclude that $\Gamma(W^{2k,p}(V))$ is reflexive and consequently, $W^{2k,p}(V)$ is also reflexive.

$(iii)$ {\it Separateness for $1\leq p<+\infty$.}

Obviously, the product space $E=L^p(V)\times \mathscr{V}^p\times\cdots\times L^p(V)\times \mathscr{V}^p\times L^p(V)$ is separable. As a consequence,
$\Gamma(W^{2k,p}(V))$ is also separable, since every closed subspace of a separable Banach space is separable. Therefore $W^{2k,p}(V)$ is separable.

The case $m=2k+1$ for some integer $k\geq 1$ can be treated similarly as the case $m=2k$ and the details are omitted.
$\hfill\Box$\\

\begin{remark}
	An important significance of Theorem \ref{wmp} is that it holds without any further restriction on the connected locally finite graph $G=(V,E)$. Not only does it not need the measure $\mu$ to have a positive lower bound, but also it does not need the weight $w_{xy}$ to satisfy symmetry. This will be very convenient for further applications, especially if one considers the $p$-Laplace equations on locally finite graphs. For such kinds of equations, Han and Shao \cite{Han-Shao} obtained existence results analogous to those of Zhang and Zhao \cite{Zhang-Zhao} under certain restrictions on the graph.
\end{remark}

Next, we shall prove Theorem \ref{wmpp} on the function space $W^{m,\mathcal{P}}(V)$, where $m\geq 1$ is an integer and $\mathcal{P}=(p_0,p_1,\cdots,p_m)$ with $p_j\geq 1$, $j=0,1,\cdots,m$. Recall that $W^{m,\mathcal{P}}(V)$ is a linear space including all functions $u$ satisfying
$$
\|u\|_{W^{m,\mathcal{P}}(V)}=\sum_{j=0}^m\|\nabla^ju\|_{p_j}<+\infty.
$$
Since the space $W^{m,\mathcal{P}}(V)$ is just $W^{m,p}(V)$, if $p_0=\cdots=p_m=p$ for some $1\leq p\leq +\infty$, we call $W^{m,\mathcal{P}}(V)$ a generalized Sobolev space. Although the computations in this space is more complicated than those in $W^{m,p}(V)$, it is still essentially based on properties of both $L^p(V)$ and $\mathscr{V}^p$.
\\

{\it Proof of Theorem \ref{wmpp}.}

$(i)$ Let $\{u_n\}$ be a Cauchy sequence in $W^{m,\mathcal{P}}(V)$, where $\mathcal{P}=(p_0,\cdots,p_m)$ with
$1\leq p_j\leq +\infty$, $j=0,1,\cdots,m$. Then $\{\nabla^ju_n\}$ is a Cauchy sequence in $L^{p_j}(V)$
if $j$ is even, or a Cauchy sequence in $\mathscr{V}^{p_j}$ if $j$ is odd. Hence there exist $u\in L^{p_0}(V)$, $f_j\in L^{p_j}(V)$ for even $j$ with $j\geq 2$ and $\mathbf{f}_j\in \mathscr{V}^{p_j}$ for odd $j$, such that $u_n\rightarrow u$ in $L^{p_0}(V)$, $\nabla^ju_n\rightarrow f_j$ in $L^{p_j}(V)$ if $j$ is even with $j\geq 2$ and $\nabla^ju_n\rightarrow \mathbf{f}_j$ in $\mathscr{V}^{p_j}$ if $j$ is odd. Clearly, $u_n\rightarrow u$ in $L^{p_0}(V)$ implies that $u_n(x)\rightarrow u(x)$ pointwisely for $x\in V$, and thus $f_j(x)=\nabla^ju(x)$ or $\mathbf{f}_j(x)=\nabla^ju(x)$ for all $j=0,1,\cdots,m$ and $x\in V$. Therefore $u\in W^{m,\mathcal{P}}(V)$ and $u_n\rightarrow u$ in $W^{m,\mathcal{P}}(V)$.

$(ii)$ Let $1<p_j<+\infty$ for all $j=0,1,\cdots, m$. Without loss of generality, we only deal with the case $m=2k$. Let $E=L^{p_0}(V)\times \mathscr{V}^{p_1}\times\cdots\times \mathscr{V}^{p_{2k-1}}\times L^{p_{2k}}(V)$. Clearly $E$ is a reflexive Banach space. Define an operator
$\Gamma:W^{2k,\mathcal{P}}(V)\rightarrow E$ by $\Gamma u=(u,\nabla u,\cdots,\Delta^ku)$. This operator satisfies that $\|\Gamma u\|_E=\sum_{j=0}^{2k}\|\nabla^ju\|_{p_j}=\|u\|_{W^{2k,\mathcal{P}}(V)}$,
and thus $\Gamma:W^{2k,\mathcal{P}}(V)\rightarrow E$ is an isometry. Since $W^{2k,\mathcal{P}}(V)$ is a Banach space,
$\Gamma(W^{2k,\mathcal{P}}(V))$ is a closed subspace of $E$. Again by Proposition 3.20 in \cite{Brezis}, we conclude $\Gamma(W^{2k,\mathcal{P}}(V))$ is reflexive and consequently, $W^{2k,\mathcal{P}}(V)$ is also reflexive.

$(iii)$ Let $1\leq p_j<+\infty$ for all $j=0,1,\cdots, m$. Without loss of generality, we still present the details of the case $m=2k$ an example. Obviously, the product space $E=L^{p_0}(V)\times \mathscr{V}^{p_1}\times L^{p_2}(V)\times\cdots\times \mathscr{V}^{p_{2k-1}}\times L^{p_{2k}}(V)$ is separable. As a consequence,
$\Gamma(W^{2k,\mathcal{P}}(V))$ is also separable, since every closed subspace of a separable Banach space is separable.
Therefore $W^{2k,\mathcal{P}}(V)$ is separable. $\hfill\Box$\\

Finally, we consider the function spaces $W_0^{m,p}(V)$ and $W_0^{m,\mathcal{P}}(V)$, which are completions of $C_c(V)$ under the norms $\|\cdot\|_{W^{m,p}(V)}$ and $\|\cdot\|_{W^{m,\mathcal{P}}(V)}$ respectively.\\

{\it Proof of Corollaries \ref{w0mp} and \ref{w0mpp}.} Obviously, $W_0^{m,p}(V)$ and $W_0^{m,\mathcal{P}}(V)$ are closed subspaces of $W^{m,p}(V)$ and $W_0^{m,\mathcal{P}}(V)$ respectively. Therefore, the properties of $W_0^{m,p}(V)$ and $W_0^{m,\mathcal{P}}(V)$ in Corollary \ref{w0mp} and \ref{w0mpp} can be deduced directly from Theorems \ref{wmp} and \ref{wmpp}.
$\hfill\Box$\\

\section{Sobolev embeddings}\label{sec5}

In this section, we shall study Sobolev inequalities on locally finite graphs.
Unlike Euclidean spaces or Riemannian manifolds, a locally finite graph has no concept of dimension. This leads us to not being able to express a function $u:V\rightarrow\mathbb{R}$, belonging to some Sobolev space, by a potential integral of its gradient, and then use H\"older's inequality to get a higher order integrability. For this reason, we add certain assumptions on measures or weights of locally finite graphs to get Sobolev inequalities in Theorems \ref{mu0} and \ref{w0}.
\\

{\it Proof of Theorem \ref{mu0}.}

For any $j_0\in\{0,1,\cdots,m\}$, if the corresponding $p_{j_0}$ satisfies $1\leq p_{j_0}<+\infty$, we have
$$
\|u\|_{W^{m,\mathcal{P}}(V)}=\sum_{j=0}^m\|\nabla^ju\|_{p_i}
\geq\mu_0^{1/{p_{j_0}}}|\nabla^{j_0}u|(x),\quad \forall x\in V.
$$
This implies
$$\|\nabla^{j_0}u\|_\infty\leq \mu_0^{-1/p_{j_0}}\|u\|_{W^{m,\mathcal{P}}(V)}.$$
On the other hand, if $p_{j_0}=+\infty$, there obviously holds
$$\|\nabla^{j_0}u\|_\infty\leq \sum_{j=0}^m\|\nabla^ju\|_{p_i}=\|u\|_{W^{m,\mathcal{P}}(V)}.$$
The above two estimates confirm \eqref{infty-mp} immediately.

Suppose that for some $j\in\{0,1,\cdots,m\}$, the corresponding $p_j$ satisfies $1\leq p_j<+\infty$. Then $\forall q\in[p_j,+\infty)$. It follows from \eqref{infty-mp} that
\begin{eqnarray*}
\|\nabla^ju\|_q&=&\le(\int_V|\nabla^ju|^qd\mu\ri)^{1/q}\\
&\leq&\|\nabla^ju\|_\infty^{1-p_j/q}\le(\int_V|\nabla^ju|^{p_j}d\mu\ri)^{1/q}\\
&\leq&C\|u\|_{W^{m,\mathcal{P}}(V)}.
\end{eqnarray*}
Thus \eqref{q-mp} is proved. $\hfill\Box$\\

{\it Proof of Theorem \ref{w0}.}

Since $W_0^{m,1}(V)$ is the completion of $C_c(V)$ under the norm $\|\cdot\|_{W^{m,1}(V)}$, it is sufficient to prove \eqref{Sob-2} for all $u\in C_c(V)$.

{\it Case 1. $m=1$.} In this case, the integer $j\in [0,1/2]$ must be $0$. For any $u\in C_c(V)$, we can choose a point $O\in V$ such that $u(O)=0$. For any given $x\in V$, there should exist a path without loop $\gamma:\{1,2,\cdots,k+1\}\rightarrow V$, such that $\gamma(i)=x_i$, $x_i\sim x_{i+1}\in E$, $x_1=x$ and $x_{k+1}=O$. It follows that
\begin{eqnarray*}
|u(x)|&=&|u(x)-u(O)|\\
&\leq& |u(x_1)-u(x_2)|+\cdots+|u(x_k)-u(x_{k+1})|\\
&\leq& \sqrt{\f{2\mu_1}{w_0}}\le\{\sqrt{\f{w_{x_1x_2}}{2\mu(x_1)}}|u(x_1)-u(x_2)|+\cdots+\sqrt{\f{w_{x_kx_{k+1}}}{2\mu(x_k)}}
|u(x_k)-u(x_{k+1})|\ri\}\\&\leq&\sqrt{\f{2\mu_1}{w_0}}\le(|\nabla u|(x_1)+\cdots+|\nabla u|(x_{k})\ri)\\
&\leq&\sqrt{\f{2\mu_1}{w_0\mu_0^2}}\le(\mu(x_1)|\nabla u|(x_1)+\cdots+\mu(x_k)|\nabla u|(x_{k})\ri)\\
&\leq&\sqrt{\f{2\mu_1}{w_0\mu_0^2}}\int_V|\nabla u|d\mu.
\end{eqnarray*}

{\it Case 2. $m>1$.} Since $u\in C_c(V)$, we have that $\Delta^ju\in C_c(V)$ for any integer $j\in [0,m/2]$. By Case 1, we can find some constant $C$
depending only on $w_0$, $\mu_0$, $\mu_1$ and $m$ such that
$$\|\Delta^ju\|_\infty\leq C\int_V|\nabla\Delta^ju|d\mu.
$$
This ends the proof of the theorem. $\hfill\Box$ \\

Before ending this section, let us make a note on the Sobolev space $W_0^{m,2}(V)$, which is obviously a Hilbert space with an inner product
$$
\langle u,v\rangle=\sum_{j=0}^m\int_V\nabla^ju\nabla^jv d\mu,\quad u,v\in W_0^{m,2}(V).
$$
As we mentioned before, the space is a completion of $C_c(V)$ under the norm $\|u\|_{W^{m,2}(V)}=\sum_{j=0}^m\|\nabla^ju\|_2$, or an equivalent norm $(\sum_{j=0}^m\|\nabla^ju\|_2^2)^{1/2}$. A natural question is whether there exists a simpler equivalent norm for this Sobolev space.
The answer to this question will be very useful when we consider related partial differential equations on locally finite graph.

Suppose that $V$ is a connected locally finite graph with symmetric weights, i.e. $w_{xy}=w_{yx}$ for all $xy\in E$. If $m=2$, obviously we have $\Delta u\in L^2(V)$, $\forall u\in W_0^{2,2}(V)$. Take a sequence $\{\phi_k\}\subset C_c(V)$ such that $\phi_k$ converges to $u$ in $W_0^{2,2}(V)$. Integration by parts gives
$$
\int_V\phi_k\Delta ud\mu=-\int_V\nabla\phi_k\nabla ud\mu.
$$
Taking the limit $k\rightarrow\infty$, we obtain
\begin{equation}\label{gai1}\int_V|\nabla u|^2d\mu\leq \f{1}{2}\le(\int_Vu^2d\mu+\int_V|\Delta u|^2d\mu\ri).\end{equation}
Hence $W_0^{2,2}(V)$ has an equivalent norm
$$\|u\|_{W_0^{2,2}(V)}=\le(\int_V(u^2+|\Delta u|^2)d\mu\ri)^{1/2}.$$
If $m=3$, there holds $\nabla\Delta u\in \mathscr{V}^2$, $\forall u\in W_0^{3,2}(V)$. Take a sequence $\{\psi_k\}\subset C_c(V)$
such that $\psi_k$ converges to $u$ in $W_0^{3,2}(V)$. Integration by parts leads to
$$-\int_V\nabla \psi_k\nabla \Delta ud\mu=\int_V\Delta \psi_k\Delta ud\mu.$$
Taking $k\rightarrow\infty$, we have
$$\int_V|\Delta u|^2d\mu\leq \f{1}{2}\le(\int_V|\nabla u|^2d\mu+\int_V|\nabla\Delta u|^2d\mu\ri),$$
which together with \eqref{gai1} gives an equivalent norm of $W_0^{3,2}(V)$
$$\|u\|_{W_0^{3,2}(V)}=\le(\int_V(u^2+|\nabla \Delta u|^2)d\mu\ri)^{1/2}.$$
By an simple induction argument, we conclude that for any integer $m\geq 1$, $W_0^{m,2}(V)$ has an equivalent norm
$$\|u\|_{W_0^{m,2}(V)}=\le(\int_V(u^2+|\nabla^mu|^2)d\mu\ri)^{1/2}.$$

\section{Further questions}\label{sec6}

As a research field with substantial content and applications, there are still many important and interesting problems about Sobolev spaces on graphs. In this section, we list two related issues that we are concerned about.\\

{\it Problem 1. For the two spaces $W_0^{m,p}(V)$ and $W^{m,p}(V)$, are they exactly the same for locally finite graphs?}

{\it Problem 2. Is there a Sobolev inequality without restrictions on $\mu$ and $w$ as in Theorems \ref{mu0} and \ref{w0}?}\\

As far as we know, the answer to the first problem is subtle. Indeed, the two spaces are the same in some specific situations, but may be different in general. Let us give some propositions and examples for these two aspects.

\begin{proposition}\label{Lp-dense}
	$L^p(V)$ is a completion of $C_c(V)$ under the norm $\|\cdot\|_p$, where $1\leq p<+\infty$ and the graph $G=(V,E)$ is an arbitrary connected locally finite graph.
\end{proposition}
{\it Proof.} Let $O\in V$ be a fixed point and $u\in L^p(V)$. Set
$$u_k(x)=\left\{\begin{array}{lll}
	u(x)&{\rm if}& x\in B_k\\[1.5ex]
	0&{\rm if}& x\in V\setminus B_k,
\end{array}\right.$$
where $B_k=\{x\in V: \rho(x)<k\}$, $k\in\mathbb{N}$.
There holds $\{u_k\}\subset L^p(V)$ and $\|u_k-u\|_p=\left(\int_{V\setminus B_k}|u|^pd\mu\right)^{1/p}=o_k(1)$, as one desired. $\hfill\Box$\\

\begin{proposition}\label{W1p-c}
	Suppose that $G=(V,E)$ is a connected locally finite graph and satisfies $(i)$ there exists a constant $\mu_0>0$ such that $\mu(x)\geq \mu_0$ for all $x\in V$;
	$(ii)$ $w_{xy}=w_{yx}$ for all $xy\in E$ and $\sum_{xy\in E}w_{xy}\leq D$ for some constant $D$ and all $x\in V$.
	Then $W^{1,p}(V)=W_0^{1,p}(V)$ for all $1\leq p<+\infty$.
\end{proposition}
{\it Proof.}
It is sufficient to show that for any $u\in W^{1,p}(V)$, there exists a function sequence $\{u_k\}\subset C_c(V)$ such that
\begin{equation}\label{tend-0}||u_{k}-u||_{W^{1,p}(V)}\rightarrow 0.\end{equation}

Fix a point $O\in V$ and let $\rho(x)$ denote the distance between $x$ and $O$. Define a sequence of cut-off functions $\eta_{k}:V\rightarrow\mathbb{R}$ as
$$\eta_{k}(x)=\left\{\begin{array}{lll} 1&{\rm if}& \rho(x)\leq k\\[1.5ex]
	0&{\rm if}&\rho(x)\geq k+1.\end{array}\right.$$
Obviously, $\{\eta_{k}\}\subset C_c(V)$, $0\leq \eta_{k}\leq 1$ and $\eta_k(x)\rightarrow 1$ for any $x\in V$ as $k\rightarrow\infty$.
For $u_{k}=\eta_{k}u$, we have $u_{k}\in C_{c}(V)$. Since $u\in L^p(V)$, there holds
\begin{equation}
	\int_{V}|u_{k}-u|^{p}d\mu=\int_{\{x:\rho(x)\geq k+1\}}|u|^pd\mu=o_k(1).\label{est-0}
\end{equation}
On the other hand, by the definition of $u_k$, there holds
\begin{eqnarray*}\nonumber
\int_{V}|\nabla u_{k}-\nabla u|^{p}d\mu&=&\underset{x\in V}\sum|\nabla u_{k}-\nabla u|^{p}(x)\mu(x)\\\nonumber
&=&\sum_{k\leq \rho(x)\leq k+1}\mu(x)|\nabla u_{k}-\nabla u|^{p}(x)+\sum_{\rho(x)>k+1}\mu(x)|\nabla u|^{p}(x)\\
&=&I_{k}+II_{k}.
\end{eqnarray*}
Since $|\nabla u|\in L^p(V)$, clearly $II_k\rightarrow 0$ as $k\rightarrow\infty$. For $I_k$, we have
\begin{eqnarray*}
I_k&=&\sum_{k\leq \rho(x)\leq k+1}\mu(x)\le|\nabla\le((\eta_k-1)u\ri)\ri|^p(x)\\
&=&\sum_{k\leq \rho(x)\leq k+1}\mu(x)\le(\sum_{y\sim x}\f{w_{xy}}{2\mu(x)}\le(\le(\eta_k(y)-1\ri)(u(y)-u(x))+
\le(\eta_k(y)-\eta_k(x)\ri)u(x)\ri)^2\ri)^{p/2}\\
&\leq&\sum_{k\leq \rho(x)\leq k+1}\mu(x)\le(\sum_{y\sim x}\f{w_{xy}}{2\mu(x)}\le(|u(y)-u(x)|+
|u(x)|\ri)^2\ri)^{p/2}\\
&\leq& 2^p\sum_{k\leq \rho(x)\leq k+1}\mu(x)\le(\sum_{y\sim x}\f{w_{xy}}{2\mu(x)}(u(y)-u(x)^2\ri)^{p/2}
+\le(\f{2D}{\mu_0}\ri)^{p/2}\sum_{k\leq \rho(x)\leq k+1}\mu(x)|u|^p(x)\\
&=&2^p\int_{\{x: k\leq\rho(x)\leq k+1\}}|\nabla u|^pd\mu+\left(\frac{2D}{\mu_0}\right)^{p/2}\int_{\{x: k\leq\rho(x)\leq k+1\}}|u|^pd\mu\\
&=&o_k(1),
\end{eqnarray*}
where we have used $|\eta_k(y)-1|\leq 1$ and $|\eta_k(y)-\eta_k(x)|\leq 1$, if $k\leq \rho(x)\leq k+1$ and $y\sim x$, and $(a+b)^2\leq 2a^2+2b^2$ for all $a,b\in\mathbb{R}$. Then we have $\lim_{k\rightarrow+\infty}\int_{V}|\nabla u_{k}-\nabla u|^{p}d\mu=0$. This together
with \eqref{est-0} leads to \eqref{tend-0}, as we desired.
$\hfill\Box$\\

Proposition \ref{W1p-c} presents a situation in which the answer of {\it Problem 1} is positive. In other situations, there are still examples that seem to support this conclusion, but are different from intuition about function spaces on Euclidean space.

Let us consider a specific locally finite graph $G=(\mathbb{N}^\ast,\mu(n),w_{n,n+1})$, where the set of vertices is $\mathbb{N}^\ast$, $\mu(n)=1/n^2$, all the neighbours of vertex $n+1$ are $\{n,n+2\}$, and $w_{n,n+1}=1$ for all $n\in\mathbb{N}^\ast$. Take a function
$$f(n)\equiv 1,\quad\forall n\in\mathbb{N}^\ast.$$
Clearly $f\in W^{1,2}(\mathbb{N}^\ast)$ with the norm $\|f\|_{W^{1,2}(\mathbb{N}^\ast)}=\sum_{n=1}^\infty 1/n^2={\rm Vol}(\mathbb{N}^\ast)$, where ${\rm Vol}(\mathbb{N}^\ast)$ denotes the volume of
$\mathbb{N}^\ast$. To find a sequence of functions that converges to $u$ in $W_0^{1,2}(\mathbb{N}^\ast)$, for any $k\in\mathbb{N}^\ast$, we define
$$f_k(n)=\le\{\begin{array}{lll}
1&{\rm if}& 1\leq n\leq k\\[1.5ex]
1-\f{j}{k}&{\rm if}& n=k+j,\,\,1\leq j\leq k\\[1.5ex]
0&{\rm if}& n>2k.
\end{array}\ri.$$
We have that $f_k$ belongs to $C_c(\mathbb{N}^\ast)$ and satisfies
$$
\int_{\mathbb{N}^\ast}|\nabla f_k|^2d\mu=\sum_{n=k}^{2k-1}(f_k(n)-f_k(n+1))^2=\f{1}{k}
$$
and
$$\int_{\mathbb{N}^\ast}(f-f_k)^2d\mu\leq 4\int_{n>k}f^2d\mu=4\sum_{n=k+1}^\infty\f{1}{n^2}.$$
It follows that
\begin{eqnarray*}
\|f-f_k\|_{W^{1,2}(\mathbb{N}^\ast)}=\|\nabla (f-f_k)\|_2+\|f-f_k\|_2
\leq\f{1}{k^{1/2}}+2\le(\sum_{n=k+1}^\infty\f{1}{n^2}\ri)^{1/2}\rightarrow 0, \quad \ {\rm as}\ \  k\rightarrow +\infty.
\end{eqnarray*}
Hence $f\in W_0^{1,2}(\mathbb{N}^\ast)$.

We can also construct a function that belongs to $W_0^{1,2}(\mathbb{N}^\ast)$ but tends to infinity as $n\rightarrow\infty$. For example, take
$$g(n)=n^{1/3},\quad\forall n\in\mathbb{N}^\ast.$$
Since
$$\int_{\mathbb{N}^\ast}|\nabla g|^2d\mu=\sum_{n=1}^\infty(g(n)-g(n+1))^2=\sum_{n=1}^\infty
\f{1}{(n^{2/3}+(n+1)^{2/3}+n^{1/3}(n+1)^{1/3})^2}$$
and
$$\int_{\mathbb{N}^\ast} g^2d\mu=\sum_{n=1}^\infty\f{1}{n^{4/3}},$$
there holds $g\in W^{1,2}(\mathbb{N}^\ast)$.
Set
$$g_k(n)=\left\{\begin{array}{lll}
n^{1/3}&{\rm if}& 1\leq n\leq k\\[1.5ex]
k^{1/3}(1-{j}/{k})&{\rm if}& n=k+j,\,\,1\leq j\leq k\\[1.5ex]
0&{\rm if}& n>2k.
\end{array}\right.$$
We have $\{g_k\}\subset C_c(\mathbb{N}^\ast)$ and a straightforward calculation gives
\begin{eqnarray*}
\int_{\mathbb{N}^\ast}|\nabla (g-g_k)|^2d\mu&\leq& 2\int_{n\geq k}|\nabla g|^2d\mu
+2\int_{n\geq k}|\nabla g_k|^2d\mu\\
&=&2\sum_{n=k}^\infty
\f{1}{(n^{2/3}+(n+1)^{2/3}+n^{1/3}(n+1)^{1/3})^2}+2\sum_{n=k}^{2k-1}\f{1}{k^{4/3}}\\
&=&O\le(\f{1}{k^{4/3}}\ri)+O\le(\f{1}{k^{1/3}}\ri)
\end{eqnarray*}
and
$$\int_{\mathbb{N}^\ast}(g-g_k)^2d\mu=\int_{n>k}(g-g_k)^2d\mu\leq \int_{n>k}g^2d\mu=
\sum_{n=k}^\infty\f{1}{n^{4/3}}=O\le(\f{1}{k^{4/3}}\ri).$$
Hence $g_k\rightarrow g$ in $W^{1,2}(\mathbb{N}^\ast)$, and thus $g\in W_0^{1,2}(\mathbb{N}^\ast)$.

Of course, in general we can construct lots of functions which are not supported in a bounded subset of a locally finite graph $G=(V,E)$ and belong to $W_0^{m,p}(V)$. But the proof of a positive answer to {\it Problem 1}  without any conditions on the graph is not easy. On the other hand, if we can construct a function in $W^{m,p}(V)$, but not in $W_0^{m,p}(V)$ for some locally finite graph $G=(V,E)$, a negative answer to {\it Problem 1} will be confirmed, which is not trivial either.

Finally, we turn to {\it Problem 2}. In Section \ref{sec5}, to prove those Sobolev inequalities, similar to Proposition \ref{W1p-c}, we need some assumptions on the measure $\mu$ and the weights $w$ of the graph. Can we get rid of or relax these constrains on $\mu$ and $w$? In fact, because of its extensive applications, it would still be of significance even if we can give any new Sobolev inequality under conditions similar to those of Theorems \ref{mu0} and \ref{w0}. This is exactly what we will do in the future. 

\section*{Acknowledgements}
This research is supported by National Natural Science Foundation of China (No. 12271039 and No. 12101355) and the Open Project Program of Key Laboratory of Mathematics and Complex System, Beijing Normal University.

\section*{Data Availability Statement}
Data sharing not applicable to this article as no datasets were generated or analysed during the current study.

\end{document}